\documentclass[%
 reprint,
 amsmath,amssymb,
 aps,
]{revtex4-1}

\usepackage{graphicx}
\usepackage{dcolumn}
\usepackage{bm}

\usepackage{amsthm}
\usepackage{color}
\newtheorem{theorem}{Theorem}

\newtheorem{lemma}{Lemma}
\newtheorem{corrolary}{Corollary}
\newtheorem{definition}{Definition}
\definecolor{mycolor}{rgb}{0.7,0.3,0.3}

\begin{document}

\preprint{APS/123-QED}

\title{Fast social-like learning of complex behaviors based on motor motifs
}

\author{Carlos Calvo Tapia$^1$, Ivan Y. Tyukin$^2$, and Valeri A. Makarov$^{1,3}$}
 	\email{Corresponding author: vmakarov@ucm.es}
	\affiliation{$^1$Instituto de Matem\'atica Interdisciplinar, Faculty of Mathematics, Universidad Complutense de Madrid, Plaza Ciencias 3, 28040 Madrid, Spain\\
$^2$University of Leicester, Department of Mathematics, University Road, LE1 7RH, United Kingdom\\
$^3$Lobachevsky State University of Nizhny Novgorod, Gagarin Ave. 23, 603950 Nizhny Novgorod, Russia}






\begin{abstract}
Social learning is widely observed in many species. Less experienced agents copy successful behaviors, exhibited by more experienced individuals. Nevertheless, the dynamical mechanisms behind this process remain largely unknown. Here we assume that a complex behavior can be decomposed into a sequence of $n$ motor motifs. Then a neural network capable of activating motor motifs in a given sequence can drive an agent. To account for $(n-1)!$ possible  sequences of motifs in a neural network, we employ the winner-less competition approach. We then consider a teacher-learner situation: one agent exhibits a complex movement, while another one aims at mimicking the teacher's behavior. Despite the huge variety of possible motif sequences  we show that the learner, equipped with the provided learning model, can rewire ``on the fly'' its synaptic couplings in no more than $(n-1)$ learning cycles and converge exponentially to the durations of the teacher's motifs. We validate the learning model on mobile robots.  Experimental results show that indeed the learner is capable of  copying the teacher's behavior composed of six motor motifs in a few learning cycles. The reported mechanism of learning is general and can be used for replicating different   functions, including, for example, sound patterns or speech.

\end{abstract}

\pacs{Valid PACS appear here}
\maketitle


\section{Introduction}

Social or imitation learning through transferring information from an experienced agent to a naive one is widely observed in many species, especially in primates and humans \cite{SocLearn}. It relies on adopting successful  actions of others and plays a fundamental role in development and communication.  Electrophysiological studies suggest that the so-called mirror neurons take place in imitation learning \cite{mirror1,mirror2}. Yet, the neuronal dynamical mechanisms underlying such a learning remain largely unknown. 

Learning a task from scratch, i.e., without any prior knowledge is a complicated problem \cite{Schaal99}. Humans rarely attempt to do it. Instead, they usually extract chains of movement primitives from instructions and demonstrations by other humans. Such a paradigm should be also implemented in future completely autonomous robots \cite{SocLearnRobot}, which will facilitate their social acceptance \cite{accept}.

While learning, an agent (human or robot) first has to isolate essential \textit{motor motifs} from the movement of another more experienced agent (e.g., ``move right'', ``turn left'', etc.). Then it must be able to compose a meaningful behavior from these motifs. Thus, the social learning is an example of a sequential solving of an inverse and then a forward problems. Given $n$ motifs one can build-up $(n-1)!$ different behaviors (e.g., for $n=10$ we get $362880$). At the first glance, such an explosive complexity may appear prohibitive for learning in small neural networks. 
 
In this work we present an approach to fast social-like learning of complex motor behaviors from an arbitrary number of motor motifs. It allows synthesizing behaviors in a neural network of a learning agent by observation of the dynamics of a teacher in no more than $(n-1)$ learning cycles. Such a linear (\textit{vs}. factorial) growth of the learning time is reasonable for many applications.

The inverse problem, i.e., the segmentation of complex movements into a sequence of primitives, is similar to the problem of discovering network motifs in neural networks \cite{Alon2007}, which relies on statistical analysis of graphs. To segment movement trajectories Stulp and colleagues \cite{Stulp09} proposed to use a clustering method with principal component analysis. However, even in common situations such an approach requires a huge number of demonstrations and may become unfeasible \cite{Colome2014}. Another approach reduces the segmentation problem to a sequential recognition of movements and comparing them to a library of motor motifs \cite{Schaal11}. Given that the library exists, this approach permits quite efficient sequencing. For the sake of simplicity, in our work we take a motif library for granted.

The forward problem can be conventionally subdivided into two tasks: 1) implementation of certain motifs (e.g., moving a hand from one point to another \cite{Carlos}), and 2) synthesis of complex behaviors from motor motifs. The first task can be approached by constructing dynamic movement primitives (DMP) \cite{Schaal2007}. A DMP is given by a dynamical system with a single or multiple attractors and several parameters that are adjusted to account for the trajectory exhibited by an experienced human. The DMP approach is based on the optimal control and in the case of relatively small number of degrees of freedom can be scaled to multiple but similar demonstrations, which produces more flexible solutions \cite{Matsubara11}.

\begin{figure*}[!ht]
\begin{center}\includegraphics[width=0.98\textwidth]{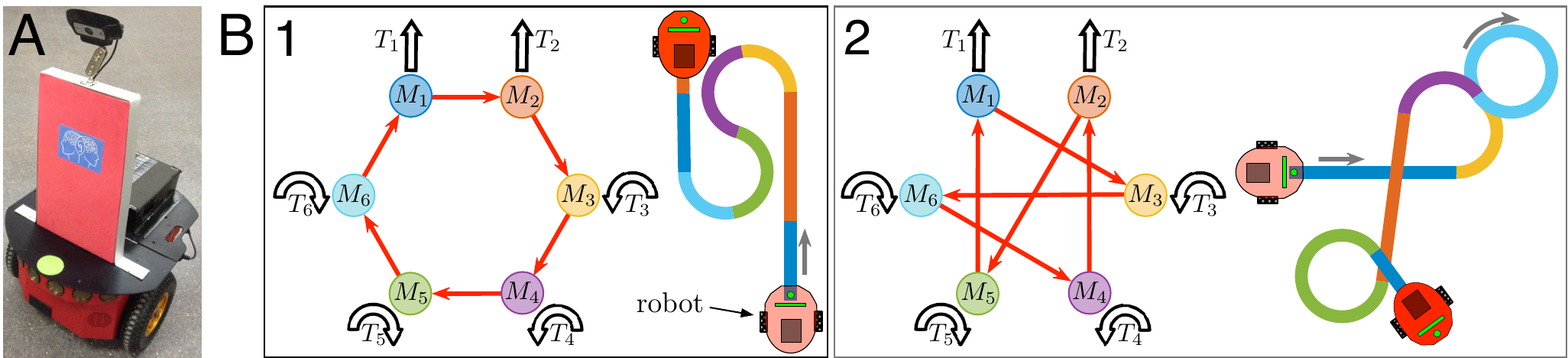}
\caption{Implementation of motor behaviors as sequences of motor motifs. A) A Pioneer 3DX mobile robot used in experiments. B) Two examples of the robot's behaviors defined by graphs of six motor motifs (see main text). Depending on the order of the motifs and their durations the robot can move along different trajectories (color codifies the active motif).}
\label{Fig_01} 
\end{center}
\end{figure*}

The second task is usually approached from the paradigm of predictive control and optimization of trajectories \cite{synth1,synth2}. We, however, explore the complementary but different question: How can a small neural network reproduce and learn one of $(n-1)!$ behaviors? From the dynamical systems viewpoint the social learning can be reduced to the problem of copying, in terms of  synchronization, ordered in time patterns of neuronal activity from one experienced agent to another, novel one. Although the general ``hardwired'' architectures of the teacher's and learner's networks can be the same, their intrinsic dynamics defining the activity patterns can differ significantly. Then the learner should be able to adjust its network couplings  by only observing the teacher dynamics, which leads to synchronization of the behaviors \cite{SM}. 

To tackle the problem here we adopt a use case  approach. As a testbed we consider the generalized Lotka-Volterra model with global inhibitory couplings among neurons \cite{Rabinovich1,SM}. Biological neural networks usually exhibit high asymmetry in the coupling patterns \cite{Benito,Divers}, i.e., reciprocal  connections between neurons $i$ and $j$ may differ significantly. Earlier it has been shown that this can be crucial for  establishing temporal associations \cite{PRL86,HighDim}. It also enables maintaining gaits in locomotion, retrieving ordered items from memory, and codification of information in population bursts \cite{Rio,Lobov}. An asymmetric Lotka-Volterra  model  can exhibit the so-called winner-less competition (WLC) behavior, which occurs in a vicinity of a stable heteroclinic loop connecting saddle equilibriums (see, e.g., \cite{Rabinovich1,C3,C5}). Then, in the presence of a weak perturbation (e.g., noise, for detail see \cite{AfrNoise} and references therein) a trajectory can wander in the phase space from one saddle to another, thus implementing a particular temporal pattern of neuronal excitation. This provides extremely rich behaviors even in small-size neural networks (see, e.g., \cite{Rabinovich1,WLC1,WLC2,Arena}).   

In this work we propose a learning algorithm capable of copying the behavior of one WLC neural network to another. Then, we consider a teacher-learner social-like situation. A teacher exhibits a complex behavior consisting of a sequence  of motor motifs (one of $(n-1)!$) with specific durations. Another agent, a learner, aims at mimicking the behavior. Despite the behavior complexity grows extremely rapidly with $n$, we show that the learning algorithm allows the learner to rewire ``on the fly'' the synaptic couplings in no more than $(n-1)$ learning cycles and converge exponentially fast to the durations of the teacher's motifs. We then validate the model on mobile robots.

\section{Motor motifs and behaviors}
\label{Sec_TESTBED}

Our first goal is to provide a neural network  architecture that would enable a robot (Fig. \ref{Fig_01}A, for details see \cite{Waves16}) to perform a series of simple motor actions or motifs. Then a complex motor behavior can be ``programmed'' as a sequence of motifs. For illustration, we select the following:
\begin{itemize}
\item[$\mathbf{M_1}$:] Go straight during time interval $T_1$.
\item[$\mathbf{M_2}$:] Go straight during time interval $T_2$.
\item[$\mathbf{M_3}$:] Turn left during time interval $T_3$. 
\item[$\mathbf{M_4}$:] Turn left during time interval $T_4$. 
\item[$\mathbf{M_5}$:] Turn right during time interval $T_5$. 
\item[$\mathbf{M_6}$:] Turn right during time interval $T_6$. 
\end{itemize}

For the sake of simplicity, in the description of the motifs we did not specify the linear robot velocity, the radii of turns, etc. These parameters will be essential in Sect. \ref{Sect_Experim}. Inclusion of similar motifs (e.g., $M_1$ and $M_2$) improves the robot flexibility and complexity of possible behaviors. Note also that the number and variety of motifs can be easily increased (with 6 motifs we can implement up to 120 behaviors).  

\begin{figure*}[hbt]
\begin{center}\includegraphics[width=0.8\textwidth]{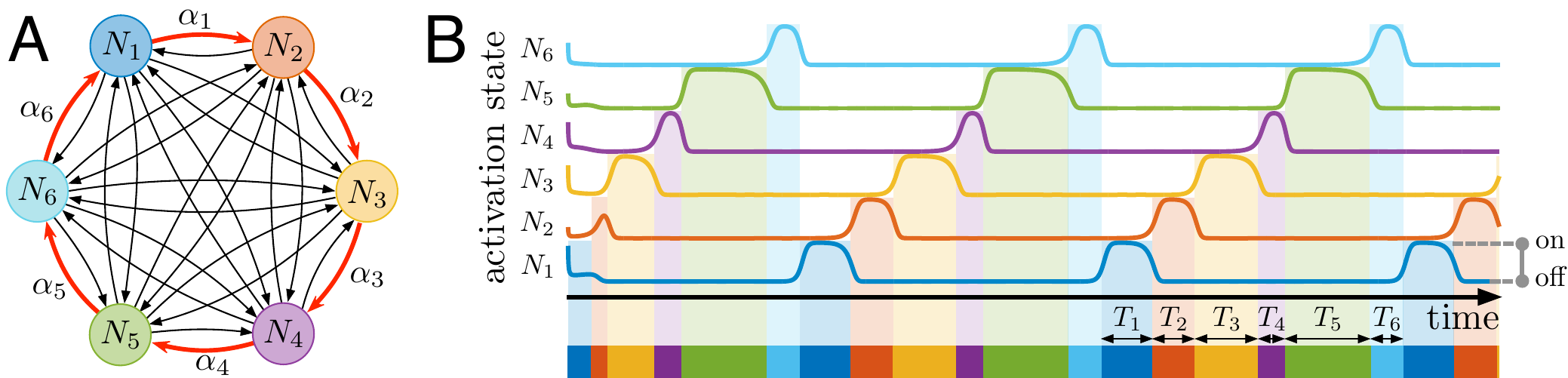}
\caption{Behavior driving neural network. A) Schematic representation of the neural network (example for $n=6$).  All neurons (circles) are mutually coupled by inhibitory synapses (black and red links). Red links mark couplings $\alpha= (0.6, 0.5, 0.7, 0.1, 0.8, 0.3)^T$ that define the activation sequence $N_1\to N_2\to N_3\to\cdots$, corresponding to the graph of motor  motifs shown in Fig. \ref{Fig_01}B.1. The coupling strength of black links is equal to 2. B) WLC dynamics generated by the network. Each neuron switches between on (active) and off (inactive) states. Bottom colored stripe shows the activation  intervals $T_1,\ldots,T_6$ of the corresponding neurons determined by  $\alpha$.}
\label{Fig02}
\end{center}
\end{figure*}

We can now connect the selected motifs in a graph and thus implement a complex behavior. Figure \ref{Fig_01}B shows two examples of graphs and the corresponding behaviors of the robot. In the first case, we selected the cyclic sequence (Fig. \ref{Fig_01}B.1): $M_1\to M_2 \to M_3 \to M_4 \to M_5 \to M_6 \to M_1 \to \cdots$, which generates a zigzag-like trajectory. In the second case the sequence was (Fig. \ref{Fig_01}B.2): $M_1\to M_3 \to M_6 \to M_4 \to M_2 \to M_5 \to M_1 \to \cdots$, which generates a more evolved trajectory.

Our second goal is to achieve a social-like learning of behaviors. We thus assign one ``more experienced'' robot as a teacher. The teacher knows how to reproduce some complex behavior given by a graph, similar to those shown in Fig. \ref{Fig_01}B. Another robot, a learner, at the beginning cannot replicate the teacher's behavior. Its neural network is initialized randomly and hence generates an arbitrary behavior (e.g., go straight or make circles). Then the learner has to learn: 1) the teacher's graph of motifs and 2) the time intervals $T_1,\ldots,T_n$ of all motifs.

\section{Neural network driving behaviors}
\label{Sect_NetModel}

Let us now provide a neural network model capable of driving a single robot, i.e., implementing all possible graphs of motor motifs $M_1,\ldots,M_n$ with the corresponding time intervals $T_1,\ldots,T_n$.
 
We consider a neural network consisting of $n\ge 3$ globally coupled neurons (Fig. \ref{Fig02}A, $n=6$). Activation of neuron $i$ during time interval $T_i$ evokes the motor motif $M_i$ and hence the robot performs the corresponding motor action while the neuron stays active. Thus, we separate the physical implementation of each motor motif from the decision-making provided by the behavior driving neural network. This enables employing the same dynamical principles for controlling different platforms and behaviors.  

To describe a sequence of motor motifs (Fig. \ref{Fig_01}B) corresponding to a sequence of activations of neurons (Fig. \ref{Fig02}A), it is convenient to introduce a directed cycle graph $G = (V,E)$ of length $n=|V|$, with the following vertices and edges (circles and red arrows in Fig. \ref{Fig_01}B):
\begin{equation}
V = \{1,\ldots,n\}, \ \ \ E = \{(i,j) \in V\times V: \ i \to j\}.
\end{equation}
Then every posible closed circuit of activation of neurons can be defined in terms of the adjacency or \textit{pathway matrix}:
\begin{equation} 
\label{Def_W}
W = \left(w_{ij}\right)_{i,j = 1}^n, \ \ \
	w_{ij} = \left\{
	\begin{array}{ll}
		1 & \mbox{ if } (j,i) \in E\\
		0 & \mbox{ otherwise. }
	\end{array}
	\right.
\end{equation}
In Eq. (\ref{Def_W}) $w_{ij}=1$ corresponds to a situation when in the  motor sequence after motif $j$ goes motif $i$. Note that $\deg^+(v) = \deg^-(v) = 1$ for all $v\in V$, therefore $W$ is an orthogonal permutation matrix.

There exist $(n-1)!$ different complete cyclic sequences represented by  hamiltonian cycles over $V$. Thus, even in a relatively small graph there can exist a huge number of possible cycles and hence motor behaviors.

To describe the network dynamics implementing sequences of motor motifs, we use the Lotka-Volterra equation generalized to an arbitrary number of neurons $n\ge3$  \cite{Rabinovich1,SM}:
\begin{equation}
\label{MainEq}
\dot{x} = x\odot(1_n - \rho_\alpha x) + \epsilon \, 1_n,  
\end{equation}
where $x(t)\in [0,1]^n$ is the activation state of the neurons resembling their firing rates, $\odot$ stands for the Hadamard product, $1_n = (1,1,\ldots,1)^T$, $\epsilon$ is a small constant ($0 < \epsilon \ll 1$), and $\rho_\alpha \in \mathcal{M}_{n\times n}(\mathbb{R}_+)$ is the interneuron coupling matrix describing global (all-to-all) synaptic links among the neurons. Moreover, since all elements of  $\rho_\alpha$ are positive, the couplings among neurons are inhibitory, i.e., activation of one neuron depresses the activity of the others.

Earlier it has been shown that if the matrix $\rho_\alpha$ satisfies appropriate conditions, then the network exhibits a winner-less competition behavior \cite{Afr03}. In particular, to satisfy the conditions we can set the connectivity matrix to (see \cite{SM} for details):
\begin{equation}
\label{Matrix_rho}
	\rho_{ij} = \left\{
	\begin{array}{ll}
	    1        & \mbox{ if } i = j      \\
		\alpha_j & \mbox{ if } w_{ij} = 1 \\
		2    & \mbox{ otherwise},
	\end{array}
	\right.
\end{equation}
where $\alpha_j \in (0,1)$. Then in Eq. (\ref{MainEq}) with $\epsilon = 0$, there exists an attracting heteroclinic circuit visiting in a cycle $n$ saddles, each of which corresponding to activation of one neuron \cite{Afr03}. Such a cycle is defined by the pathway matrix $W$ and corresponds to the graph $G$ (Fig. \ref{Fig_01}B). Since the heteroclinic circuit is structurally unstable, any small perturbation ($\epsilon >0$) leads to an emergence of a stable limit cycle in the vicinity of the destroyed heteroclinic loop, which is the only attractor implementing a WLC behavior. 

Figure \ref{Fig02}B illustrates an example of the WLC dynamics. For the sake of visual clarity we have implemented the simplest excitation circuit: $N_1\to N_2 \to N_3 \to N_4 \to N_5\to N_6$ corresponding to the graph shown in Fig. \ref{Fig_01}B.1.  Loosely speaking, at each time instant only one neuron is on and it inhibits the activity of other neurons until another neuron ``switches'' on and the process repeats. The duration of the on-state of neuron $j$ is an increasing function of $\alpha_j$,  which tends to infinity if $\alpha_j \to 1$ \cite{SM}. Thus, we can individually control the time intervals of all motifs.

Therefore, playing with the pathway matrix $W$ we can select one of the $(n-1)!$ sequences of activation of the neurons (and hence the graph of motor motifs), while choosing $\alpha = (\alpha_1,\ldots,\alpha_n)^T$ we can set the temporal extension of  activation for each neuron (and hence the time intervals $T_1,\ldots,T_n$ for all motifs). Thus, the model (\ref{MainEq}) satisfies our needs and offers an extremely flexible way for generating a variety of motor behaviors.

\section{Fast learning of behaviors in a teacher-learner framework}
\label{Sect_FastLearn}

Let us now consider a social-like situation with two agents: an experienced robot (a teacher) and a naive one (a learner). The teacher has a well-trained neural network (\ref{MainEq}), which allows reproducing some complex behavior (see, e.g., Fig. \ref{Fig_01}B). As we have shown in Sect. \ref{Sect_NetModel} the teacher's behavior is defined by the pathway matrix $W_x$ (subindex $x$ refers to the teacher) and the vector $\alpha$. 

The dynamics of the neural network driving the learner is given by:
\begin{equation}
\label{EqLearner}
	\dot{y} = y\odot(1_n - \rho_\gamma y) + \epsilon \, 1_n, 
\end{equation}
where $y(t)\in [0,1]^n$ is the activation state of the neurons and $\rho_\gamma(t)$ is the coupling matrix with entries $\gamma \in \mathbb{R}_+^n$ [counterpart of $\alpha$ in Eqs. (\ref{MainEq}) and (\ref{Matrix_rho})].

At the beginning the learner has an arbitrary pathway matrix $W_y$ (counterpart of $W_x$) and the time interval vector $\gamma$. Thus, the dynamics of the learner can differ significantly from the dynamics of the teacher. Then the learner's goal is to learn the teacher's behavior and reproduce it. 

\subsection{Learning time intervals}
\label{Sect_LearnTime}

Let us first assume that the learner ``knows'' the pathway matrix $W_x$, i.e., the graph of motor motifs of the teacher $G$. However, the time intervals of the motifs ($T_1,\ldots,T_n$) given by $\alpha$ in the teacher and by $\gamma(t)$ in the learner are unknown. Thus, we have: 
\begin{equation}
W_y(t) = W_x, \ \ \ \gamma(0) \ne \alpha.
\end{equation}
Note that the initial values $\gamma(0)$ are  taken from a random distribution and we can even have a situation with $\gamma_j(0) \ge 1$ for some $j\in V$. In such a case the WLC conditions are not satisfied and the network evolves to an equilibrium state, and hence the learner either stays still or performs a single motor motif.     
  
The learner aims at replicating the teacher's behavior, i.e., $y(t)$ for $t>t'>0$ must repeat in some sense the behavior of $x(t)$. To proceed further let us introduce two essential definitions.
   
\begin{definition} 
\label{Def1}
Under the learning of a temporal pattern we understand a situation when a learner, sharing the same pathway matrix with the teacher ($W_y=W_x$), independently on the initial conditions $x(0)$, $y(0)$, $\gamma(0)$, can tune its synaptic couplings $\gamma(t;x)$ in such a way that
\begin{equation}
\label{learn_def}
	\lim_{t \to \infty} \gamma(t; x) = \alpha. 
\end{equation}   
\end{definition}

We note that condition (\ref{learn_def}) implies that the learner will eventually succeed in reproducing the temporal sequence shown by the teacher:
\begin{equation}
\label{expo_def}
\lim_{t \to \infty} (x(t) - y(t-t^*)) = 0_n,
\end{equation}   
where $t^*$ is some phase lag. This happens because for a fixed (stationary) $\gamma \in (0,1)^n$ the WLC model (\ref{EqLearner}) has a single structurally stable attractor. Thus, given (\ref{learn_def}) at $t\to\infty$ the learner and the teacher will have the same stable limit cycle in their phase spaces. Then, the phase lag $t^*$ in (\ref{expo_def}) corresponds to a delay in movement along the limit cycle of the teacher and of the learner.

\begin{definition} 
\label{Def2}
We say that the learning of a temporal pattern is exponentially fast if
\begin{equation}\nonumber
	\|\gamma(t) - \alpha \|_2 \le M e^{-\kappa t} , \ \ \forall t \ge 0
\end{equation} 
for some constants $M>0$ and $\kappa>0$.
\end{definition}  
We note that exponential convergence as opposed to mere asymptotic one ensures robustness of the learning process.

\subsubsection{Learning rule and exponential convergence}
Following the adaptive system approach \cite{Tyukin1,Tyukin2}  we now  postulate the learning rule:
\begin{equation} 
\label{Lrule}
	\gamma(t) = \gamma_0 + W_y^T\big(\theta(t) - (x(t) - x_0)\big),
\end{equation}
where $\theta(t)\in \mathbb{R}^n$ is the control variable obeying:
\begin{equation} 
\label{Def_th}
\dot{\theta} = x \odot(1_n - \rho_{\gamma} x) + \epsilon \, 1_n, \ \ \ \theta(0)=0.
\end{equation}
Note that learning rule (\ref{Lrule}), (\ref{Def_th}) does not depend on the direct knowledge of the teacher's internal parameters $\alpha$ and even of the own behavior of the learner $y(t)$, but it is fully based on the observation of the teacher's dynamics  $x(t)$. Thus, the learning can be ``mental'', i.e., without motor execution.  

For further calculations it is convenient to introduce the quadratic term:
\begin{equation}
\label{Def_P}
p(t) = W_y^Tx(t)\odot x(t).
\end{equation}
Then we can formulate the following result.

\begin{theorem} 
\label{Th1}
Under learning rule  (\ref{Lrule}), (\ref{Def_th}) the learner learns exponentially fast any temporal pattern exhibited by  the teacher.

Moreover, assuming that the teacher exhibits a stable $T$-periodic pattern (i.e., $x(t+T) = x(t)$), we have the following  estimate for the convergence exponent: 
\begin{equation}
\label{LearnExp}
\kappa = \min \left \{ \langle p_j \rangle \right \}, \ \ \langle p \rangle = \frac{1}{T}\int_{t_0}^{t_0+T} p(t)\, \mathrm{d}t.
\end{equation}

\end{theorem}

The proof is provided in Appendix \ref{app:1}.

\subsubsection{Numerical simulations}
Let us now illustrate Theorem \ref{Th1} and the learning abilities of the proposed learning rule. For the sake of visual clarity we built 3-neuron networks for the teacher and the learner. We then assigned three motor motifs: $M_1$ - move right, $M_2$ - turn left, and $M_3$ - turn right. The teacher's sequence of motor motifs was set by selecting the pathway matrix:
\begin{equation}
W_x = \left(
\begin{array}{ccc}
0 & 1 & 0\\
0 & 0 & 1\\
1 & 0 & 0
\end{array}
\right),
\end{equation} 
which imposes the motif sequence $M_1 \to M_3 \to M_2\to \cdots$. The time interval vector was set to: $\alpha = (0.2, 0.6, 0.8)^T$, i.e., $T_1<T_2<T_3$. Then, we set the initial values of the learner's couplings $\gamma(0) = (1.6, 0.1, 2.3)^T$. Note that in this case $\gamma_{1}, \gamma_3>1$ and hence the learner cannot exhibit a WLC dynamics at the beginning.  We now numerically integrate the model (\ref{MainEq}), (\ref{EqLearner}), (\ref{Lrule}), (\ref{Def_th}) and by using $x(t)$ estimate the convergence exponent (\ref{LearnExp}): $\kappa = 0.0134$. 

\begin{figure}[!t]
\centerline{
\includegraphics[width=0.44\textwidth]{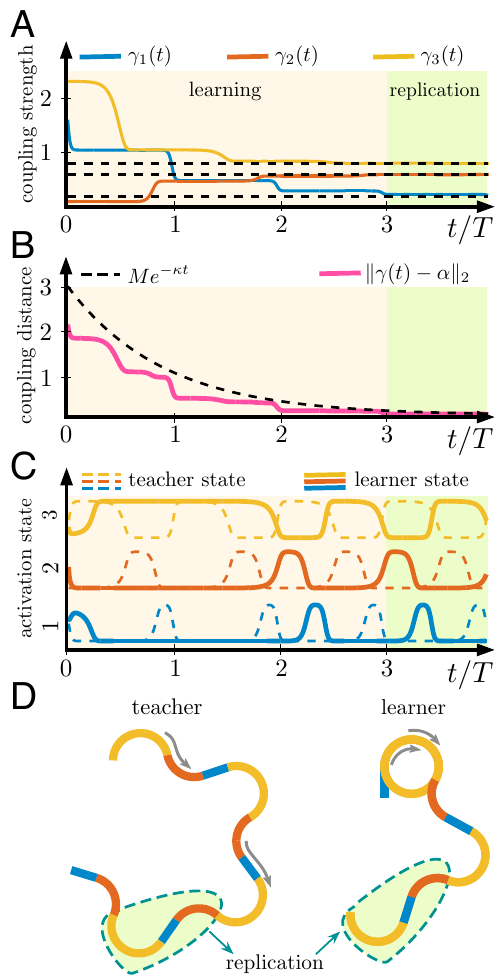}
}
\caption{Representative example of exponentially fast learning of temporal patterns. A) Dynamics of the coupling strengths of the learner $\gamma(t)$ (solid curves) converging to the corresponding couplings of the teacher $\alpha$ (dashed lines). B) Exponentially fast convergence of the couplings. C) Synchronization (with a time  shift) of oscillations in the learner with the teacher. D) Simulation of the robots' movements. The teacher (left) performs zigzag movements, while the learner (right) at the beginning moves in a circle, but after two learning cycle it starts replicating the teacher's behavior (color coding is the same as in C).}
\label{Fig_03}
\end{figure}

Figure \ref{Fig_03} illustrates the learning process. As expected the couplings $\gamma(t)$ converge exponentially to $\alpha$ (Figs. \ref{Fig_03}A, \ref{Fig_03}B). Step-like changes occur when the corresponding component of $p(t)$ goes through its maximum. This happens in time intervals when the activation passes from one neuron to another. At the beginning the learner does not follow a WLC dynamics, as expected, neuron 3 (responsible for right turn) stays active for a long time interval (Fig. \ref{Fig_03}C, $0.3T \lessapprox  t \lessapprox 1.8T$, where $T$ is the oscillation period).  Nevertheless, once all $\gamma_j$ enter the region $(0,1)$ (Fig. \ref{Fig_03}A, $t\approx 1.5T$) the learner starts generating a WLC behavior that finally synchronizes with the teacher's dynamics (Fig. \ref{Fig_03}C, $t\gtrapprox 2T$). 

Figure  \ref{Fig_03}D shows simulations of the robots' movements driven by the neural networks (see also Fig. \ref{Fig_01}B). The teacher performs a zigzag movements. Since the time interval $T_3$ corresponding to motif $M_3$ (turn right) is the longest, the robot exhibits a tendency to turn right (Fig.  \ref{Fig_03}D, teacher). At the beginning the learner has neuron 3 activated (Fig. \ref{Fig_03}C, $0.3T \lessapprox  t \lessapprox 1.8T$) and hence it performs circular movements (Fig.  \ref{Fig_03}D, learner). Then, however, it synchronizes with the teacher and replicates its behavior.

\subsection{Learning a sequence of motor motifs}

Let us now consider the general case. At the beginning the learner starts from an arbitrary orthogonal matrix $W_y(0)\in \mathcal{M}_{n\times n}(\{0,1\})$, such that $W_y(0)\ne W_x$. Then, $W_y(0)$ induces a wiring in the graph of motor motifs  different from that of the teacher, i.e., $E_y \ne E_x$. 

Without loss of generality we can assume that the teacher exhibits a $T$-periodic dynamics: 
\begin{equation}
x(t)=x(t+T),
\end{equation}
with a given pathway matrix $W_x$ and a coupling vector $\alpha$. Following the derivation similar to the proof of Theorem \ref{Th1} (see Appendix \ref{app:1}) we obtain the following equation describing the dynamics of the learner's  couplings:

\begin{equation}
\label{Gam_dyn2}
\dot{\gamma} = -p\odot\gamma + q,
\end{equation}
where $p$ is defined by (\ref{Def_P}) and 
\begin{equation}
\label{EQ:19b}
	q = W_y^Tx \odot \left[2x - (2 - W_y^TW_x\alpha)\odot W_y^TW_xx \right].
\end{equation}
Then, from (\ref{Gam_dyn2}) we get
\begin{equation}
\label{General_Solution}
\gamma(t) = \left(\gamma_0 + g(t)\right)\odot f(t),
\end{equation}
where 
\begin{equation}
\label{f_and_g}
f(t) = e^{-\int_0^t p(\tau)\, \mathrm{d}\tau}, \ \ \ g(t) = \int_0^t{q(\tau)\oslash f(\tau)\, \mathrm{d}\tau},
\end{equation}
where $\oslash$ stands for the Hadamard division. We note that $f(t)$ is a strictly positive function and $\gamma(t)$ does not now converge to $\alpha$, but oscillates due to $W_y\ne W_x$.

\subsubsection{Criterium of successful learning}

For a successful learning (Def. 1), the pointwise limit of $\gamma(t+kT)$ [see also (\ref{Def_H})] for $k\to\infty$ must be equal to $\alpha$, regardless  the value of $t$. However, in general it is an oscillating function given by (\ref{General_Solution}). We then introduce the following mean squared deviations $e_j:\mathbb{R}^2\to \mathbb{R}$:
\begin{equation}\label{Def_e}
	e_j(\delta) = \frac{1}{2}\Big\langle\big(\delta_{1} - ( \gamma_j - \delta_{2}f_j )\big)^{2}\Big\rangle, \ \ \forall j \in V,
\end{equation}
where $\delta = (\delta_1,\delta_2)^T$ and $\langle \cdot \rangle$ denotes the time averaging operator over the period $T$. Note that for a given $\delta$ the learner can evaluate (\ref{Def_e}) along the trajectory, i.e., it is an observable variable.   

We now can find the minimum of $e_j$ by evaluating the gradient $\nabla e_j$ and the Hessian matrix:
\begin{equation}\label{GRHES}
\mathcal{H}(e_j) =  \left(\hspace{-0.1cm}
\begin{array}{cc}
	1 					& \langle f_j \rangle \\
	\langle f_j \rangle & \langle f_j^2 \rangle
\end{array}\hspace{-0.1cm}\right).
\end{equation}
We observe that $\det(\mathcal{H}) = \mathrm{Var}[f_j] > 0$. Thus, $e_j$ is convex and reaches the global minimum at $\nabla e_j=0$. Solving this equation, we get the coordinates of the minimum:
\begin{equation}
\label{Linear_System}
\delta_j^* = \mathcal{H}^{-1}\left( \hspace{-0.1cm}
\begin{array}{c}
\langle \gamma_j \rangle\\ 
\langle \gamma_jf_j \rangle 
\end{array}
\hspace{-0.1cm}\right) = (\langle \gamma_j \rangle - \beta_j\langle f_j  \rangle , \beta_j)^T,
\end{equation}
where $\beta_j = \mathrm{Cov}[\gamma_j,f_j]/\mathrm{Var}[f_j]$. In other words, $\delta_j^*$ is composed of the intercept and slope values of linear regression of $\gamma_j$ over $f_j$. 

Given this observation we can formulate the following result.

\begin{theorem}
\label{THEOREM2}
Let neuron $i\in V$ be the successor of neuron $j$ in the learner's sequence of motor motifs, i.e. $(j,i)\in E_y$. Then
\begin{equation}
\label{EqCritConv}
(j,i) \in E_x \Leftrightarrow  e_j(\delta_j^*) = 0.
\end{equation}
\end{theorem}
The proof is given in Appendix \ref{app:3}. Equation (\ref{EqCritConv}) provides a criterium that enables indirect inferring on the structure of the pathway matrix  in the teacher $W_x$ by the learner by evaluating the error $e(\delta_j^*)$.

\subsubsection{Iterative adjustment of $W_y$}

Based on Theorem \ref{THEOREM2}, we propose an iterative scheme to adjust the pathway matrix of the learner $W_y$. For further calculations, it is convenient to introduce the discrete time $k$, such that the time-line is divided into intervals:
\begin{equation}
I_k=[(k-1)T,kT), \ \ \ k = 1,2,\ldots
\end{equation} 
We also denote:
\begin{equation}
x_{k} = x(kT), \  \theta_{k} = \theta(kT), \  \gamma_{k} = \gamma(kT).
\end{equation}

The learner starts with a random orthogonal matrix $W_y[0]$ and the time interval vector $\gamma_0$. Then in each interval $I_k$  we have the matrix $W_y[k-1]$ that defines the structure of the coupling matrix $\rho_{k-1}(\gamma)$.  The dynamics of the learner  [compare to (5), (9), and (10)] is given by ($t\in I_k$):    
\begin{equation} \label{k-system}
	\left\{
	\begin{array}{l}
		\dot{y} = y\odot(1 - \rho_{k-1}(\gamma)y) + \epsilon\, 1_n \\
		\dot{\theta} = x\odot(1 - \rho_{k-1}(\gamma)x) + \epsilon \, 1_n \\
\gamma = \gamma_{k-1} + W_y^T[k-1]\big(\theta - \theta_{k-1} - x + x_{k-1}\big).
	\end{array}
	\right.
\end{equation}
We now find $\delta^*_{j}$ using (\ref{Linear_System}) and then evaluate $e_j(\delta_j^*)$ from (\ref{Def_e}). By applying Theorem \ref{THEOREM2} (i.e., detecting $j$ such that $e_j(\delta_j^*)\ne 0$), we find edges in the graph $G_y$ that should  be rewired (they do not appear in $G_x$). After rewiring we get new matrix $W_y[k]$. This procedure is repeated for $k=1,2,\ldots$. Thereby, we obtain a sequence of the matrices $\{W_y[k]\}$. 

Below we will propose an algorithm of the network rewiring that ensures:
\begin{equation}
W_y[k] = W_x, \ \  \ \forall k \ge k^* >0,
\end{equation}  
where $k^* < n$, i.e., the number of steps required to find $W_x$ grows no more than linearly with the number of motor motifs $n$. Thus, the learner can learn rapidly  the pathway matrix of the teacher. Then, in accordance with Theorem \ref{Th1}, the learner also learns the time interval vector $\alpha$ exponentially fast.

\subsection{Learning algorithm}
\label{LearALG}
Summarizing the abovementioned results we  propose the following algorithm of learning.

\begin{figure*}[!ht]
\centerline{
\includegraphics[width=0.8\textwidth]{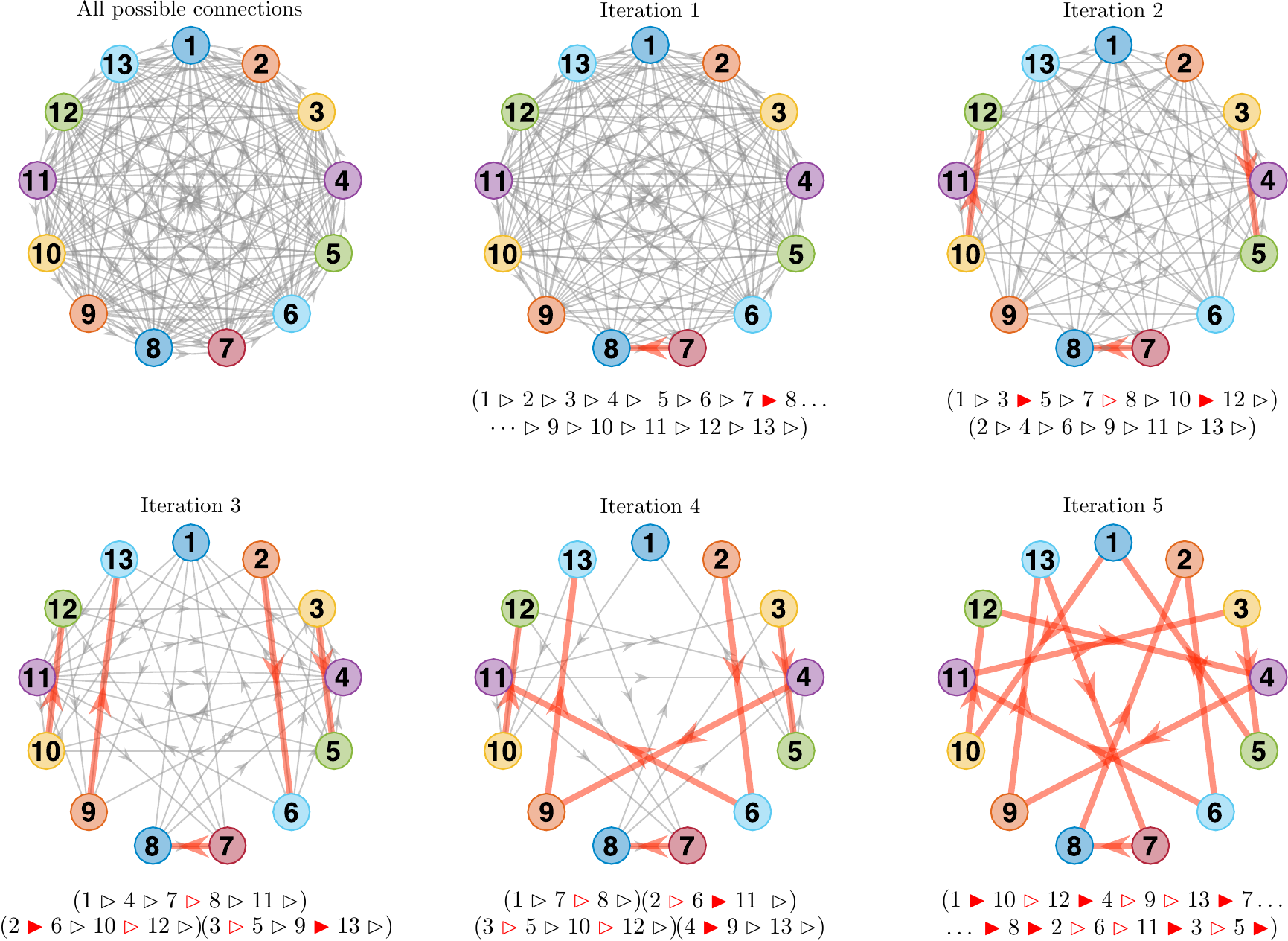}
}
\caption{Representative example of dynamical learning of the teacher's graph $1\to 10 \to 12 \to 4 \to 9 \to 13 \to 7 \to 8 \to 2 \to 6 \to 11 \to 3 \to 5 \to \cdots$ defining the sequence of motor motifs in the teacher's behavior. First panel shows in gray all possible connections. Iterations from 1 to 5 show correctly identified connections in red and the remaining ones in gray. A legend below each panel shows the learner's graph. Filled red triangles correspond to newly identified connections, open red triangles mark previously identified connections, and open black traingles label wrong connections.}
\label{fig_LERN}
\end{figure*}

\begin{itemize}
\item[$\mathbf{S_0}$:] Initialize the parameters:
1) Select the initial permutation $\sigma : V \to V$, $\sigma = (1,2,\ldots, n)$; 2) Select random $\gamma_0$; 3) $\theta_0 \leftarrow 0$; 4) $k \leftarrow 0$; 5) $\Omega \leftarrow V$.

\item[$\mathbf{S_1}$:] While $\Omega \neq \emptyset$ repeat:
\begin{enumerate}
	\item $k \leftarrow k + 1$.
	\item Set $W_y[k-1] = \left(w_{ij}\right)_{i,j=1}^n$, where
	\begin{equation}\label{Updated_W}
		w_{ij} = \left\{
		\begin{array}{ll}
   			 1    & \mbox{ if } \ i = \sigma(j)      \\
				0    & \mbox{ otherwise.}
		\end{array}
		\right.
	\end{equation}
	\item Integrate (\ref{k-system}) over $I_k$ by using the coupling matrix $\rho_{k-1}$ given by $W_y[k-1]$.
	\item Find $\delta_j^*$ and $e_j(\delta_j^*)$ from (\ref{Def_e})--(\ref{Linear_System}).
	\item Set $\Omega = \{j \in V: e_j(\delta_j^*) \neq 0\}$ and label its elements in ascending order $\Omega= \{\omega_1,\dots,\omega_m\}$, with $\omega_k < \omega_{k+1}$, $ \forall k < m$.
	\item For each $i \in \{1,\dots,m\}$ update the permutation $\sigma(w_i) \leftarrow \sigma(w_{i^+})$, where
	\begin{equation}
	\label{iPermut}
			i^+ = \left\{\begin{array}{ll}
				i+1  & \mbox{if } \ i < m \\
				1    & \mbox{if } \ i = m.
			\end{array}\right.
	\end{equation}
\end{enumerate}

\item[$\mathbf{S_2}$:] Using the last $W_y$ (note, that now $W_y = W_x$) apply the learning rule (\ref{Lrule}), (\ref{Def_th}) for $t\ge kT$.
\end{itemize}

The learning time of the proposed algorithm is given by the following theorem.

\begin{theorem}
\label{Theorem3}
The iterative learning process described at step $\mathbf{S}_1$ of the algorithm converges in $(n-1)$ iterations at most.
\end{theorem}
The proof is provided in Appendix \ref{ProofTheorem3}. Thus, the learning time grows linearly with the number of neurons, while the number of possible behaviors grows extremely fast, as factorial.

Let us now illustrate the algorithm numerically. We built two neural networks with 13 neurons each (a teacher and a learner) and selected a random pathway matrix $W_x$ for the teacher. Then, we  initialized the learner (step $\mathbf{S}_0$)  and started the algorithm.  Figure \ref{fig_LERN} illustrates the process of learning of the connectivity matrix $W_y$ (step $\mathbf{S}_1$). 

At the beginning there exist $4.79\times 10^8$ possible graphs of the motor motifs (panel ``all possible connections''). At  the first iteration (lasting one period of the teacher's behavior) the learning algorithm identifies that the link going from neuron 7 to neuron 8 exists in the teacher, while the others defined by $W_y[0]$ do not (Fig. \ref{fig_LERN}, Iteration 1). Thus, we keep link $7\to8$ and change the others ($\mathbf{S}_1$, item 6 of the algorithm). Then, a new wiring is set up  and we repeat the calculation over the second period of the teacher's behavior. Within Iteration 2 the learner finds two more correct links: $10\to 12$ and $3\to 5$ (Fig. \ref{fig_LERN}, Iteration 2). In the following three iterations, the learner gets to learn all connections $1\to 10\to 12\to \cdots$. Thus, it was able to learn the teacher's sequence out of $4.79\times 10^8$ possible cases in only five iterations or demonstrations of the teacher's behavior.

\begin{figure*}[hbt]
\begin{center}\includegraphics[width=0.65\textwidth]{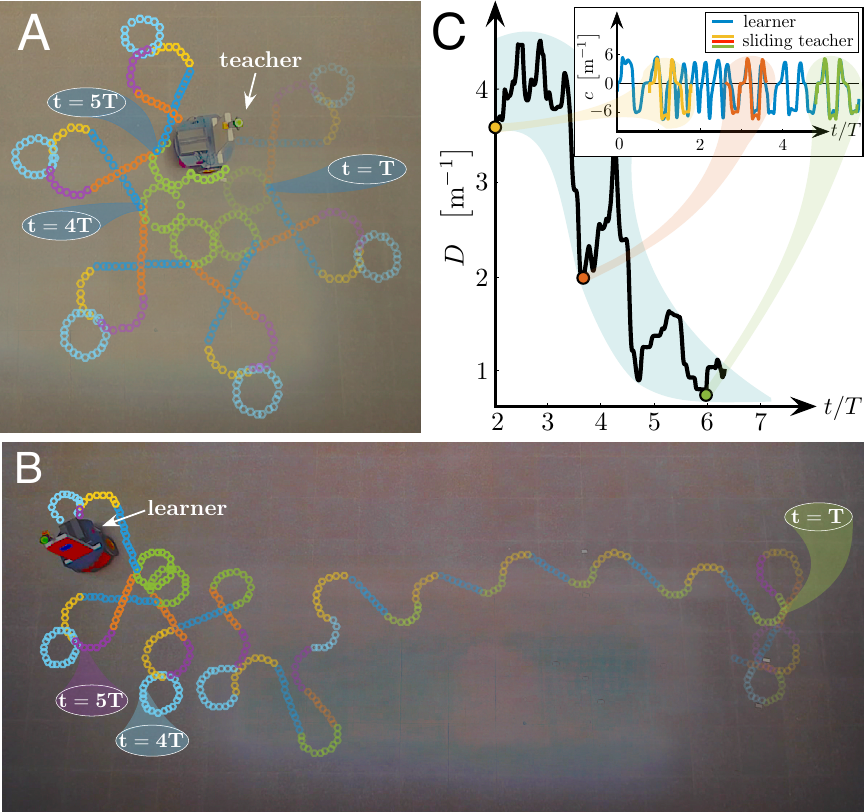}
\caption{Experimental validation of learning. A) Robot-teacher implements the behavior composed of six motor motifs shown in Fig. \ref{Fig_01}B.2. The brightness of colors corresponds to time (the brighter, the closer to the present). B) Trajectory of the robot-learner. At the beginning it differs significant from the teacher's behavior. However, in only three cycles the robot learns the sequence of turns and then two more cycles are needed  for final adjustment of the coupling strengths. Eventually the learner replicates almost exactly the behavior of the teacher. C) Quantification of the learning error (\ref{Metric}), defined as the distance between trajectories of the teacher and learner. Inset illustrates the trajectory curvature of the learner (blue curve) with superimposed curvatures of the teacher in three time windows marked by circles of different colors. The level of coincidence increases with time.}
\label{Fig04}
\end{center}
\end{figure*}

Then in step $\mathbf{S}_2$ the time intervals of the motor motifs (the values of $\gamma$) converge exponentially (Theorem \ref{Th1}). Note that this process does not depend on the number of motifs. Thus, the time interval required to learn a complex motor behavior grows linearly with the complexity of the behavior (i.e., with the number of motifs).

\section{Experimental validation}
\label{Sect_Experim}

Let us now test experimentally the proposed learning model. We equipped two robotic platforms (Pioneer 3DX, Adept Mobilerobotics,  linear sizes l$\times$w: $45.5\times 38.1$ cm, Fig. \ref{Fig_01}A) with neural networks consisting of $n=6$ neurons. The network dynamics had been implemented in Matlab (Mathworks) running on an on-board Intel NUC PC connected to the platform through a COM interface. One platform played the role of a teacher (with a predesigned complex motor behavior), while the other one was designated as a learner with the driving network randomly initialized. 

In both platforms we implemented six motor motifs provided in Sect. \ref{Sec_TESTBED}. The linear robot velocity was set to $10$ cm/s and the radii of turns to $17$ cm. These parameters were chosen  to fit the robot's behavior within the available arena (8$\times$3.5 m). During all experiments the robots' behavior was recorded by a wide-angle zenithal USB camera. The geometric distortions of the actual position of the robots in the arena (due to projection) were corrected by applying an appropriate transformation.     

The neural network dynamics of the teacher has been translated into motor commands by 
\begin{equation}
i(t) = \arg \hspace{-0.4mm} \max \{x_j(t)\}_{j=1}^6,
\end{equation}
where $i(t)\in\{1,2,\ldots,6\}$ is the number of the motor motif performed by the robot at time instant $t$. The same mechanism has been used to drive the learner.  

The pathway matrix $W_x$, defining the coupling matrix $\rho_\alpha$, and hence the sequence of the motor motifs of the teacher, was selected such that it reproduced the graph shown in Fig. \ref{Fig_01}B.2 (i.e., $M_1 \to M_3 \to M_6 \to M_4 \to M_2 \to M_5$). The values of $\alpha$ were set in such a way that the durations of the motor motifs were $T_1=7.0$ s, $T_2=7.1$ s, $T_3=4.1$ s, $T_4=4.1$ s, $T_5=9.4$ s, $T_6=11.0$ s. Then, during one period $T=\sum T_i =42.7$ s the teacher exhibited a behavior similar to that shown schematically  in Fig. \ref{Fig_01}B.2.  Figure \ref{Fig04}A shows the teacher's trajectory corresponding to five repetitions of the motor pattern. The robot repeats the behavior shown in Fig. \ref{Fig_01}B.2 and produces a complex ``flower''--like trajectory (Fig. \ref{Fig04}A).

The neural network of the learner was initialized at random  ($W_y$ and $\gamma$ were arbitrary chosen). Then the learning algorithm (Sect. \ref{LearALG}) has been activated and the robot started moving. Figure \ref{Fig04}B shows the robot trajectory. At the beginning the robot performed quite chaotic movements, only partially reproducing the teacher's motor pattern. However, after three periods the learner was able to capture completely the pathway matrix $W_x$ (i.e., $W_y(t)=W_x$ for $t > 4T$) and started repeating the teacher pattern quite reliably. Then, the values of $\gamma$ were also tuned (exponentially fast) and finally the learner was able to reproduce the teacher's behavior faithfully.

To quantify the difference between the trajectories performed by the robots  we introduced the following metric:
\begin{equation}
\label{Metric}
D(t) = \min_{\tau \in [0,T]}\left [\frac{1}{T}\int_{t-T}^t|c_T(s)-c_L(s-\tau)|^2\, \mathrm{d}s\right ]^\frac{1}{2}, 
\end{equation}
where $c_{T}$ and $c_{L}$ are the curvatures of the teacher's and learner's trajectories, respectively, evaluated by:
\begin{equation}
\label{Curvat}
c(t) = \frac{\tilde{x}'(t)\tilde{y}''(t)-\tilde{y}'(t)\tilde{x}''(t)}{(\tilde{x}'(t)^2+\tilde{y}'(t)^2)^{3/2}},
\end{equation}  
where $(\tilde{x}(t),\tilde{y}(t))$ is the parametric robot trajectory captured by the video camera. Since the linear robot velocity is constant, curvature (\ref{Curvat}) describes uniquely the robot trajectory on a 2D plane. We also note that the metric (\ref{Metric}) is invariant with respect to the translational and rotational symmetries, and the phase lag in the learner [see Eq. (\ref{expo_def})]. Thus, the introduced metric fulfills all requirements for a distance measure between two trajectories. 
      
Experimentally acquired  trajectories are not precise due to a number of reasons, starting form errors in the odometry \cite{Waves16}, identification of the robot position, and ending by the distortions produced by the camera. Thus, to evaluate the metric (\ref{Metric}) we upsampled the trajectories shown in Figs. \ref{Fig04}A and \ref{Fig04}B and denoised the data by a continuous wavelet approach \cite{WavDen}. Figure \ref{Fig04}C shows the dynamics of the distance between the trajectories of the teacher and the learner. As expected the distance strongly oscillates at the beginning due to rewiring and changes of $\gamma$, but then  decreases and approaches a stationary level, defined by measurement and motor noise. Thus, the learner robot can indeed adapt ``on the fly'' the wiring of its decision-making network and mimic the movements of the teacher.

\section{Conclusion}

Cognitive learning of motor behaviors is a complicated problem even for humans. In this work, we have proposed an efficient mechanism for dynamical learning in neural networks. In a social-like situation, a teacher (mobile robot in experiments) can exhibit a complex behavior consisting of a sequence of $n$ motor motifs (simple motor actions). In general, there exist $(n-1)!$  combinations of motifs. Moreover, all motifs have specific time durations that also play an important role in formation of the final motor behavior. Thus, the complexity and variety of the motor patterns available to an agent grow extremely rapidly with $n$. This poses two problems: 1) How can a small neural network implement such a variety of different behaviors? and 2) How can a learner (i.e., another  neural network) copy the unique behavior of the teacher starting from random initial conditions?  

To implement a motor behavior (e.g., to drive a robot-teacher) we have proposed to link the output of $n$ neurons with $n$ motor motifs. Then, an activation of a given neuron evokes the execution of a corresponding motor motif. This approach allows separating the ``behavioral'' neural network from the executive motor  part. Thus, the network dynamics does not depend on specific motor motifs and physical properties of the agent, and hence can be easily transferred to other technical devices. To avoid contradicting commands in such an architecture, two neurons cannot be active simultaneously. This property has been achieved by employing the winner-less competition paradigm that ensures sequential ``switching'' of neurons in the network and also enables simple independent control of the durations of ``on'' states of the neurons \cite{SM,C3,Rabinovich1}. 
 
Using the motor motif approach, a learner has to solve two tasks: 1) Find the teacher's graph of motifs (i.e., the pathway matrix), and 2) Adjust the time durations of motifs (i.e., the coupling strengths). To deal with the latter task, earlier we provided a neural network model that enabled learning of the activation times \cite{SM,MC}. However, to reproduce the  teacher's behavior it was necessary to start simulation from appropriate initial conditions (although rather general). Besides, learning required a long time. Here, we have introduced a novel learning rule and proved its exponentially fast (and thus robust) convergence independently on the initial conditions (Theorem 1). We note that the learning rule does not depend on the learner's state variable. Thus, during the learning an explicit motor implementation of the behavior is not necessary. Therefore, the agent can simply ``observe'' the teacher without moving itself. Such a silent or ``mental'' learning can be observed in humans and also in artificial cognitive agents \cite{IEEE13}.   

The learning of the teacher's graph (task 1) is based on the dynamic evaluation of the fitness of interneuron couplings over the learner's trajectory. We have provided a rigorous approach (Theorem 2) that enables network rewiring ``on the fly''. Such a rewiring requires at most $(n-1)$ iterations of the algorithm for its complete convergence (Theorem 3). Thus, the learning time grows only linearly with the number of motor motifs. In a numerical simulation with networks of $n=13$ neurons we observed that the learner could find the teacher's pathway matrix (one out of $4.79\times 10^8$ possible) in five steps only. 

Thus, the provided approach enables fast learning in social-like situations. To validate it, we tested the algorithm on mobile robots. First, we implemented six motor motifs (go straight, turn left, etc.) as executable commands. These motifs can be considered as ``genetically'' programed primitives available to a roving robot. Then, we implemented neural networks in two robots and designated one of them as a teacher and the other as a learner. The teacher has been programed to reproduce a complex flower-like trajectory, while the learner was set at random. Then, the learner started to learn the teacher's trajectory. At the beginning its movements were quite chaotic. However, after few cycles the learner has successfully ``copied'' the teacher's behavior. To confirm quantitatively the learning behavior, we have proposed a special metric, invariant with respect to translational and rotational symmetries, and to a time lag in the execution of trajectories. We have shown that the distance between the teacher's trajectory and the learner's behavior decreases in time.          

In conclusion, the reported mechanism of learning is quite general and can be used for replication of different behaviors on different platforms. Moreover, the  behaviors to be replicated need not to be motor. For instance, one can also think about replication of sound patterns or speech. It may also serve as a linker, connecting different scenarios and behaviors during cognitive navigation \cite{BC15}.


\begin{acknowledgments}
This work has been supported by the Russian Science Foundation under project 15-12-10018 (the problem statement and theoretical development) and by the Spanish Ministry of Economy and Competitiveness under grant FIS2014-57090-P.
\end{acknowledgments}

\appendix

\section{Proof of Theorem \ref{Th1}} 
\label{app:1}
First let us recall that the common operation of multiplication of matrices has a higher priority then the Hadamard operations. Now deriving (\ref{Lrule}) and using Eqs. (\ref{MainEq}), (\ref{Def_th}) we obtain:  
\begin{equation}
\label{General_case}
\dot{\gamma} = -W^T\big( x\odot (\rho_\gamma - \rho_\alpha)x\big).
\end{equation}
Then we observe that $(\rho_\gamma - \rho_\alpha)x =  W\big((\gamma - \alpha)\odot x\big)$. Substituting this to (\ref{General_case}) we get:
\begin{equation}
\label{Eq_inter1}
\dot{\gamma} = -W^Tx \odot (\gamma - \alpha) \odot x. 
\end{equation}
Therefore,
\begin{equation}
\label{Imposition_gam}
\frac{d(\gamma-\alpha)}{dt} = - p \odot (\gamma - \alpha),
\end{equation}
where $p\in \mathcal{C}\big(\mathbb{R},[0,1]^n\big)$, given by (\ref{Def_P}), is a function of time completely defined by the dynamics of the teacher. Integrating (\ref{Imposition_gam}) we obtain:
\begin{equation} 
\label{EXP_CONV}
\gamma(t) - \alpha  = (\gamma_0 - \alpha) \odot e^{-\int_0^t{p(\tau)d\tau}}.
\end{equation}
Since $x(t) = 0$ is not a solution of (\ref{MainEq}) and hence $x(t)>0$ under a WLC dynamics, there exists $\kappa >0$ such that $p(t) \ge \kappa$, which ensures the exponential convergence: $\|\gamma(t) - \alpha\|  \le \|\gamma_0 - \alpha \| e^{-\kappa t}$. If $x(t)$ is periodic, then a lower bound for $\kappa$ is given by (\ref{LearnExp})$_\square$

\section{Proof of Theorem \ref{THEOREM2}}
\label{app:3}

\subsection{Preliminary results}
For further calculations it is convenient to introduce two constants:
\begin{equation}
A = f(T) \in (0,1)^n, \ \ \ B = g(T)\in \mathbb{R}^n,
\end{equation} 
where $f$, $g$ are given by (\ref{f_and_g}), and $T$ is the oscillation period of the teacher's pattern. These constants $A$ and $B$ are completely determined by the teacher's dynamics and the pathway matrix of the learner, $W_y$.

\begin{lemma}
\label{LM1}
Let $T$ be the period of the teacher's pattern $x(t+T)=x(t)$. Then, under the learning rule (\ref{Gam_dyn2}),(\ref{EQ:19b}) we have
\begin{equation}
\label{LemmaEq}
\gamma(t + kT) = \gamma(t) - \big(1_n - A^{\circ k}\big)\odot \big(\gamma_0 - A\odot B\oslash (1_n - A)\big)\odot f(t),
\end{equation} 
\end{lemma}
Note that $(\cdot)^{\circ k}$ in (\ref{LemmaEq}) stands for the Hadamard $k$-th power.

\

\noindent
\textit{Proof:}
We notice that $p(t)$ [see Eq. (\ref{Def_P})] is $T$-periodic and thus
\begin{equation*}
\int_0^{t+kT}p(\tau)\,\mathrm{d}\tau =  k\int_0^Tp(\tau)\,\mathrm{d}\tau+\int_0^tp(\tau)\,\mathrm{d}\tau.
\end{equation*}
Therefore,
\begin{equation}\label{L1}
	f(t+kT) = A^{\circ k}\odot f(t).
\end{equation}

Then, $q(t)$ [see Eq. (\ref{EQ:19b})] is also $T$-periodic, which together with (\ref{L1}) yields 
\begin{multline}\label{L2}
	g(t+kT) = \sum_{l = 0}^{k-1}\int_{lT}^{(l+1)T}q(\tau)\oslash f(\tau)d\tau +\\
 {\displaystyle + \int_{kT}^{kT+t}q(\tau)\oslash f(\tau)d\tau = }\\
	= {\displaystyle B\odot\sum_{l = 0}^{k-1}A^{\circ -l} +  A^{\circ -k} \odot \int_0^t{q(\tau)\oslash f(\tau)d\tau}}.
\end{multline}

Now by substituting (\ref{L1}) and (\ref{L2}) into  (\ref{General_Solution}) we get
\begin{multline}
	\gamma(t+kT) =\gamma_0 \odot A^{\circ k}\odot f(t) + \\
	 +{\displaystyle B\odot\sum_{l = 1}^k A^{\circ l}\odot f(t) + \gamma(t) - \gamma_0\odot f(t)} = \gamma(t)-\\
  {\displaystyle - \Bigg((1_n - A^{\circ k})\odot\gamma_0 - B\odot\sum_{l = 1}^k A^{\circ l}\Bigg)\odot f(t)}.
\end{multline}

Finally, since $A \in (0,1)^n$ we can evaluate the geometric series:
\begin{equation}
	\sum_{l = 1}^k A^{\circ l} = A \odot (1_n - A^{\circ k}) \oslash (1_n - A),
\end{equation}
which yields (\ref{LemmaEq})$_\square$

\

\begin{corrolary}
\label{CORR1} 
Given the sequence $$\left\{\gamma(t+kT) : 0 \le t < T \right\}_{k = 1}^\infty \subset \mathcal{C}\left([0,T),\mathbb{R}^n\right),$$ there exists a constant vector $C \in\mathbb{R}^n$ such that
\begin{equation}\label{Def_H}
	 \lim_{k\to\infty}\gamma(t+kT) = \gamma(t) - C\odot f(t) \hspace{0.5cm} \forall t \in [0,T)
\end{equation}
\end{corrolary}

\

\noindent 
\textit{Proof:}
Since $A \in (0,1)^n$, $\lim_{k\to\infty} A^{\circ k} = 0$. Then denoting $C = \gamma_0 - A \odot B \oslash (1_n - A)$ and applying Lemma \ref{LM1} we arrive to (\ref{Def_H})$_\square$

\subsection{Proof of the theorem}
On the one hand, given that $(j,i) \in E_x\cap E_y $ and  taking into account Eq. (\ref{EXP_CONV}) we get
\begin{equation}
\gamma_j(t+kT) = \alpha_j + \left(\gamma_j(0) - \alpha_j\right)f_j(t+kT)
\end{equation}
and thus
\begin{equation}
\lim_{k\to\infty}\gamma_j(t+kT) = \alpha_j.
\end{equation}
Then by Corollary \ref{CORR1} there exists $C_j\in\mathbb{R}$ such that
\begin{equation}
\gamma_j(t) = {\displaystyle \lim_{k\to\infty}\gamma_j(t+kT) + C_jf_j(t)} = \alpha_j + C_jf_j(t).
\end{equation}
Thus, by choosing $\delta^*_j = (\alpha_j,C_j)^T$ we obtain $e_j(\delta^*_j) = 0$. This $\delta^*_j$ matches the obtained from (\ref{Linear_System}).

On the other hand, if $e_j(\delta^*_j) = 0$ then
\begin{equation}
\gamma_j(t) = \delta_{j1}^* + \delta_{j2}^*f_j(t)\  \Rightarrow \ \dot\gamma_j(t) = \delta_{j2}\dot{f}_j(t). 
\end{equation}
Now, we notice that $\delta_{j2}^*\dot{f}_j = -p_j\delta_{j2}^*f_j = -p_j(\gamma_j-\delta_{j1}^*)$ and thus
\begin{equation}
\dot{\gamma}_j(t) = -p_j(t)\gamma_j(t) + q_j(t),
\end{equation}
where $q_j(t) = p_j(t)\delta_{j1}^*$. Let $l\in V$ be the predecessor of $i\in V$, i.e., $(l,i)\in E_x$. From Eq. (16) we get:
\begin{equation}
q_j(t) = 2p_j(t) + (\alpha_l - 2)x_i(t)x_l(t).
\end{equation}
Now comparing the expressions for $q_j$ we get 
\begin{equation}
x_l(t) =\frac{\delta_{j1}^* -2}{\alpha_l - 2} x_j(t),
\end{equation}
which yields $l=j$, $\delta_{j1}^* = \alpha_j$ and thus the $i$-th rows of the pathway matrices $W_x$ and $W_y$ are the same$_\square$

\section{Proof of Theorem \ref{Theorem3}}
\label{ProofTheorem3}

Let us define, for each vertex $j \in V$, the cumulative set of tested successors up to step $k > 0$:
\begin{equation}
S_k(j) = \cup_{l = 1}^k \{\sigma_l(j)\}.
\end{equation}
We note that $\sigma(j)$ changes (item 6 of the algorithm) ensuring an exhaustive sequential search due to the rule $\sigma(w_i) \leftarrow \sigma(w_{i+1})$ and the ascending order established for $\Omega$ (item 5 of the algorithm). If at step $k$ the tested successor matches the teacher's network (item 4 of the algorithm), then the $j$-th node is excluded from $\Omega$ at step $k+1$. Therefore, $S_{k+1}(j)$ stops increasing. Thus, we have
\begin{equation}
\label{C2}
j\in \Omega_k \Leftrightarrow | S_{k+1}(j) | = k+1.
\end{equation}
Now we notice that $S_k(j) \subset V \setminus \{j\}$ and hence the cardinality of the cumulative sets satisfies the inequality:
\begin{equation}
\label{C3}
|S_k(j)| \le |V\setminus \{j\}| = n-1, \ \ \ k = 1,2,\ldots
\end{equation}

Using (\ref{C2}) and (\ref{C3}) we thus obtain that  $\forall j\in V:$ $j \in \Omega_k$ iff $k < n -1$. Therefore, $\Omega_{n-1} = \emptyset$, i.e., the iterative learning process converges at most in $n - 1$ steps$_\square$


\begin{thebibliography}{99}
\bibitem{SocLearn} 
W. Hoppitt and K.N. Laland. \textit{Social Learning: An Introduction to Mechanisms, Methods, and Models} (Princeton Univ. Press, 2013).

\bibitem{mirror1}
V. Gallesea and A. Goldmanb. Mirror neurons and the simulation theory of mind-reading. Trend. Cogn. Sci. \textbf{4}, 252 (2000).


\bibitem{mirror2} R. Cook, G. Bird, C. Catmur, C. Press, and C. Heyes. Mirror neurons: From origin to function. Behav. Brain Sci. \textbf{37}, 177 (2014).

\bibitem{Schaal99} 
S. Schaal. Is imitation learning the route to humanoid robots? Trend. Cogn. Sci. \textbf{3}, 233 (1999).


\bibitem{SocLearnRobot}
C. Breazeal, D. Buchsbaum, J. Gray, D. Gatenby, and B. Blumberg. Learning from and about others: Towards using imitation to bootstrap the social understanding of others by robots. Artif. Life \textbf{11}, 31 (2005).


\bibitem{accept} E. Broadbent, R. Stafford, and B. MacDonald. Acceptance of healthcare robots for the older population: Review and future directions. Int. J. Soc. Robot. \textbf{1}, 319 (2009).


\bibitem{Alon2007}
U. Alon. Network motifs: theory and experimental approaches. Nat. Rev. Genet. \textbf{8}, 450 (2007).

\bibitem{Stulp09}
F. Stulp, E. Oztop, P. Pastor, M. Beetz, and S. Schaal.
Compact models of motor primitive variations for predictable reaching and obstacle avoidance. Proc. Conf. IEEE-RAS on Humanoid robots (2009).


\bibitem{Colome2014} 
A. Colome and C. Torras. Dimensionality reduction and motion coordination in learning trajectories with dynamic movement primitives. Proc. IEEE/RSJ Int. Conf. Intell. Robot. Syst. (2014).

\bibitem{Schaal11} 
F. Meier, E. Theodorou, F. Stulp, and S. Schaal. Movement segmentation using a primitive library. Proc. IEEE/RSJ Int. Conf. Intell. Robot. Syst., 3407 (2011).

\bibitem{Carlos}
J.A. Villacorta-Atienza, C. Calvo, S. Lobov, and V.A. Makarov. Limb movement in dynamic situations based on generalized cognitive maps. Math. Model. Nat. Phenom. \textbf{12}, 15 (2017).

\bibitem{Schaal2007} 
S. Schaal, P. Mohajerian, and A. Ijspeert. Dynamics systems vs. optimal control - a unifying view. Progr. Brain Res. \textbf{165}, 425 (2007).

\bibitem{Matsubara11}
T. Matsubara, S. Hyon, and J. Morimoto. Learning parametric dynamic movement primitives from multiple demonstrations. Neur. Netw. \textbf{24}, 493 (2011).

\bibitem{synth1} 
Y. Tassa, T. Erez, and E. Todorov. Synthesis and stabilization of complex behaviors through online trajectory optimization. IEEE/RSJ Int. Conf. Intell Robots Syst., 4906 (2012).

\bibitem{synth2}
A. Ansari and T. Murphey. Sequential action control: Closed-form optimal control for nonlinear and nonsmooth systems.  IEEE Trans. Robot. \textbf{32}, 1196 (2016).

\bibitem{SM}
A. Selskii and V. Makarov. Synchronization of heteroclinic circuits through learning in coupled neural networks. Regular and Chaotic Dynamics \textbf{21}, 97 (2016).

\bibitem{Rabinovich1}
M. Rabinovich, A. Volkovskii, P. Lecanda, R. Huerta, H. Abarbanel, and G. Laurent. Dynamical encoding by networks of competing neuron groups: Winnerless competition. Phys. Rev. Lett. \textbf{87}, 068102 (2001).

\bibitem{Benito}
N. Benito, G. Martin-Vazquez, J. Makarova, V. Makarov, and O. Herreras. The right hippocampus leads the bilateral integration of gamma-parsed lateralized information. Elife \textbf{5}, e16658 (2016).

\bibitem{Divers}
G. Martin-Vazquez, N. Benito, V.A. Makarov, O. Herreras, and J. Makarova. Diversity of LFPs activated in different target regions by a common CA3 input. Cerebral Cortex \textbf{26}, 4082 (2016).

\bibitem{PRL86}
H. Sompolinsky and I. Kanter. Temporal association in asymmetric neural networks. Phys. Rev. Lett. \textbf{57}, 2861 (1986).

\bibitem{HighDim}
I. Tyukin, A. Gorban, C. Calvo, J. Makarova, and V.A. Makarov. High-dimensional brain. A tool for encoding and rapid learning of memories by single neurons. Bull. Math. Biol. DOI: 10.1007/s11538-018-0415-5 (2018).

\bibitem{Rio}
E. Del Rio, V.A. Makarov, M.G. Velarde, and W. Ebeling. Mode transitions and wave propagation in a driven-dissipative Toda-Rayleigh ring. Phys. Rev. E \textbf{67}, 056208 (2003).

\bibitem{Lobov}
S.A. Lobov, M.O. Zhuravlev, V.A. Makarov, and V.B. Kazantsev. Noise enhanced signaling in STDP driven spiking-neuron network. Math. Model. Nat. Phenom. \textbf{12}, 109 (2017).

\bibitem{C3}
M.A. Cohen and S. Grossberg. Absolute stability of global pattern formation and parallel memory storage by competitive neural networks. IEEE Tran. Syst. Man Cyber. \textbf{13}, 815 (1983).


\bibitem{C5}
P. Varona, M. Rabinovich, A. Selverston, and Y. Arshavsky. Winnerless competition between sensory neurons generates chaos: A possible mechanism for molluscan hunting behavior. Chaos \textbf{12}, 672 (2002).

\bibitem{AfrNoise}
V.S. Afraimovich, V.P. Zhigulin, and M.I. Rabinovich. On the origin of reproducible sequential activity in neural circuits. Chaos \textbf{14}, 1123 (2004).

\bibitem{WLC1} 
F. Hadaeghi, M. Reza, H. Golpayegani, and G. Murray. Towards a complex system understanding of bipolar disorder: A map based model of a complex winnerless competition. J. Theor. Biol. \textbf{376}, 74 (2015).

\bibitem{WLC2}
T. Rost, M. Deger, and M. Nawrot. Winnerless competition in clustered balanced networks: inhibitory assemblies do the trick. Biol. Cybern. DOI: 10.1007/s0042 (2017).

\bibitem{Arena}
P. Arena, L. Fortuna, D. Lombardo, L. Patane, and M.G. Velarde. The winnerless competition paradigm in cellular nonlinear networks: Models and applications. Int. J. Circ. Theor. Appl. \textbf{37}, 505 (2009).

\bibitem{Waves16}
C. Calvo, J. Villacorta-Atienza, V.I. Mironov, V. Gallego, and V.A. Makarov. Waves in isotropic totalistic cellular automata: Application to real-time robot navigation. Advances in Complex Systems \textbf{19}, 1650012 (2016). 

\bibitem{Afr03}
V. Afraimovich, M. Rabinovich, and P. Varona. Heteroclinic Contours in Neural Ensembles and the Winnerless Competition Principle. Int. J. Bifurc. Chaos \textbf{14}, 1195 (2004).

\bibitem{Tyukin1}
I.Yu. Tyukin, D.V. Prokhorov, and C. van Leeuwen. Adaptation and parameter estimation in systems with unstable target dynamics and nonlinear parametrization. IEEE Trans.  Autom. Contr. \textbf{52}, 1543 (2007).

\bibitem{Tyukin2}
I. Tyukin. \textit{Adaptation in Dynamical Systems} (Cambridge Univ. Press, 2011).

\bibitem{WavDen}
A.E. Hramov, A.A. Koronovskii, V.A. Makarov, A.N. Pavlov, and E. Sitnikova. \textit{Wavelets in Neuroscience} (Springer, 2015).


\bibitem{MC}
V.A. Makarov, C. Calvo, V. Gallego, and A. Selskii. Synchronization of heteroclinic circuits through learning in chains of neural motifs. IFAC-PapersOnLine \textbf{49}, 80 (2016).

\bibitem{IEEE13}
J. Villacorta-Atienza and V.A. Makarov.  Neural network architecture for cognitive navigation in dynamic environments. IEEE Trans. Neur. Netw. Learn. Syst. \textbf{24}, 2075 (2013).

\bibitem{BC15}
J. Villacorta-Atienza, C. Calvo, and V.A. Makarov. Prediction-for-CompAction: Navigation in social environments using generalized cognitive maps. Biol. Cybern. \textbf{109}, 307 (2015).



\end{thebibliography}
\end{document}